\documentclass[11pt,letterpaper]{amsart}

\usepackage{amsmath,amsthm,amssymb,latexsym,amscd,graphs}
\usepackage{epic, eepic}

\newtheorem{Theorem}{Theorem}[section]
\newtheorem{Corollary}[Theorem]{Corollary}
\newtheorem{Main}[Theorem]{Main Theorem}

\newtheorem{Lemma}[Theorem]{Lemma}
\newtheorem{Proposition}[Theorem]{Proposition}

\theoremstyle{definition}
\newtheorem{Definition}[Theorem]{Definition}

\newtheorem{Question}[Theorem]{Question}

\theoremstyle{remark}

\title{Coloring Graphs with Crossings}
\author[Bogdan Oporowski]{Bogdan Oporowski$^*$}
\thanks{*Department of Mathematics, Louisiana State University, Baton Rouge, LA  70803}
\thanks{~\,{\ttfamily bogdan@math.lsu.edu}}
\author[David Zhao]{David Zhao$^\dag$}
\thanks{$\dag$Department of Computer Science, University of Texas, Austin, TX  78712}
\thanks{~\,{\ttfamily wzhao@cs.utexas.edu}}
\date{\today}


\begin{document}

\setlength{\parskip}{1.8ex}
\setlength{\parindent}{18pt}

\begin{abstract}
We generalize the Five Color Theorem by showing that it extends to graphs with two
crossings.
Furthermore, we show that if a graph has three crossings, but does not contain
$K_6$ as a subgraph, then it is also $5$-colorable.
We also consider the question of whether the result can be extended to graphs
with more crossings.
\end{abstract}

\maketitle

%%%%%%%%%%%%%%%%%%%%%%%%%%%%%%%%%%%%%%%%%%%%%%%%%%%%%%%%%%%%%%%%%%%%%%%%%%%%%%%%%%%%%%%%

\section{Introduction}

In this paper, $n$ will denote the number of vertices, and $m$ the number of edges,
of a graph $G$.
A coloring of $G$ is understood to be a {\em proper coloring}; that is, one in which adjacent vertices always receive distinct colors.

We will consider \textit{drawings\/}\index{good drawing} of graphs in the plane
${\mathbb R}^2$ for which no three edges have a common crossing.
A crossing of two edges $e$ and $f$ is {\em trivial\/} if $e$ and $f$
are adjacent or equal, and it is {\em non-trivial\/} otherwise.
A drawing is {\em good \/} if it has no trivial crossings.
The following is a well-known easy lemma.

\begin{Lemma}\label{trivial_x}
  A drawing of a graph can be modified to eliminate all of its trivial crossings, with the number of non-trivial crossings remaining the same.
\end{Lemma}

To avoid complicating the notation, we will use the same symbol for a graph and its
drawing in the plane.
We will refer to the {\em regions\/} of a drawing of a graph $G$ as the maximal open sets
$U$ of ${\mathbb R}^2-G$ such that for every two points $x,y\in U$, there exists a polygonal $xy$-curve in $U$.

\begin{Definition}
  The {\em crossing number\/} of a graph $G$, denoted by $\nu(G)$, is the minimum
  number of crossings in a drawing of $G$.
  An {\em optimal drawing\/} of $G$ is a drawing of $G$ with exactly $\nu(G)$ crossings. 
\end{Definition}

\begin{figure}
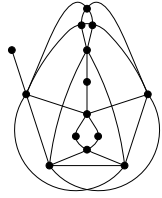

  \label{im-fig}
  \begin{graph}(2,2.45)(0,.25)
    \filledareasfalse
    \graphnodesize{.1}
    \graphlinewidth{.01}
    \allinethickness{.011cm}
    \roundnode{1}(.191,1.412)
    \roundnode{2}(1,2)
    \roundnode{3}(1.809,1.412)
    \roundnode{4}(1.5,.461)
    \roundnode{5}(.5,.461)
    \roundnode{v}(1,1.149)
    
    \roundnode{a1}(.85,.85)
    \roundnode{a2}(1.15,.85)
    \roundnode{a3}(1,.675)
    \roundnode{a4}(1,2.55)
    \roundnode{a5}(1,1.575)
    \roundnode{a6}(.93,2.33)
    \roundnode{a7}(1.07,2.33)
    \roundnode{a8}(0,2)
    
    \spline(1,1.149)(1.25,.75)(.5,.461)
    \spline(1,1.149)(.75,.75)(1.5,.461)
    \spline(1,2)(1.125,2.5)(.75,2.75)(.191,1.412)
    \spline(1,2)(.875,2.5)(1.25,2.75)(1.809,1.412)
    \area(.191,1.412){(.75,2.25,.05)(1.25,2.25,.2)(1.809,1.412,.05)}
    \edge{3}{4}
    \edge{4}{5}
    \edge{5}{1}
    \edge{v}{1}
    \edge{v}{2}
    \edge{v}{3}
    \bow{2}{4}{.05}
    \bow{2}{5}{-.05}
    \bow{1}{4}{-.5}
    \bow{5}{3}{-.5}
    \edge{a8}{1}
  \end{graph}
  \caption{A graph with an essential immersion of $K_6$}
\end{figure}

\begin{Definition}
  Suppose $G'$ and $G$ are graphs.
  A function $\alpha$ with domain $V(G') \cup E(G')$ is an {\em immersion\/}
  of $G'$ into $G$ if the following hold:
  \begin{enumerate}
  \item the restriction of $\alpha$ to $V(G')$ is an injection into $V(G)$;
  \item for an edge $e$ of $G'$ incident to $u$ and $v$, the image $\alpha(e)$ is a path
    in $G$ with ends $\alpha(u)$ and $\alpha(v)$; and
  \item for distinct edges $e$ and $f$ of $G'$, their images $\alpha(e)$ and $\alpha(f)$
    are edge-disjoint.
  \end{enumerate}
  The immersion $\alpha$ is {\em essential\/} if additionally $\alpha(e)$ and
  $\alpha(f)$ are vertex-disjoint whenever $e$ and $f$ are not adjacent, and
  it is an {\em embedding\/} if $\alpha(e)$ and $\alpha(f)$ are internally
  vertex-disjoint for all distinct $e$ and $f$.
  If $v$ is a vertex of $G$, and $\alpha$ is an essential immersion of $G'$ into $G$ such that $v=\alpha(u)$ for some vertex $u$ of $G'$, and $\alpha(e)$ is a single-edge path for each $e$ incident with $u$, then $\alpha$ is called a {\em $v$-immersion\/} of $G'$ into $G$.  
  We will also say that $\alpha$ is an immersion of $G'$ {\em onto\/} $G$
  if the range of $\alpha$ is $V(G) \cup E(G)$.
  Depending on the properties of $\alpha$, we will say that $G'$ is {\em immersed},
  {\em essentially immersed}, {\em embedded}, or {\em $v$-immersed into\/} or
  {\em onto\/} $G$.
  An example appears in Figure \ref{im-fig}.
\end{Definition}

It is worth noting that if, for every edge $e$ of $G'$, the path $\alpha(e)$
consists of a single edge, then $G'$ is a subgraph of $G$.
All immersions considered in the remainder of this paper will be essential.

\begin{Proposition}\label{num_crossings}
  If $n \ge 3$, then $\nu(G) \ge m - 3n + 6$.
\end{Proposition}

\begin{proof}
  Since $m \le 3n-6$ in a planar graph, every edge in excess of this bound introduces at least one additional crossing.
\end{proof}

\begin{Corollary}\label{k6}
  The crossing number of the complete graph $K_6$ is three.
\end{Corollary}

\begin{proof}
  It is easy to draw $K_6$ with exactly three crossings, while Proposition~\ref{num_crossings} implies that
  $\nu(K_6) \ge 3$.
\end{proof}

%%%%%%%%%%%%%%%%%%%%%%%%%%%%%%%%%%%%%%%%%%%%%%%%%%%%%%%%%%%%%%%%%%%%%%%%%%%%%%%%%%%%%%%%

\section{Immersions and Crossings}

In this section we present several results that relate crossings of a drawing with immersions of a graph.

\begin{Lemma}\label{immersion-crossing}
  Suppose $G$ is a good drawing with exactly $k$ crossings and there is 
  an essential immersion of $G'$ onto $G$.
  Then $G'$ has a good drawing with exactly $k$ crossings.
\end{Lemma}

\begin{proof}
  Let $\alpha$ be an essential immersion of $G'$ onto $G$.
  Draw $G'$ by placing each vertex $v$ at $\alpha(v)$, drawing each edge $e$ so that
  it follows $\alpha(e)$, and then perturbing the edges slightly so that no edge
  contains a vertex and no three edges cross at the same point.
  Each crossing of edges $e$ and $f$ in $G'$ arises from the corresponding paths
  $\alpha(e)$ and $\alpha(f)$ either crossing or sharing a vertex.
  In the latter case, the crossing is trivial as the immersion $\alpha$ is essential.
  The conclusion now follows immediately from Lemma~\ref{trivial_x}.
\end{proof}

Thus we have the following:

\begin{Corollary}\label{onto-greater}
  If $G'$ is essentially immersed into $G$, then $\nu(G') \le \nu(G)$.
\end{Corollary}

We may also use essential immersions to extend the Five Color Theorem.

\begin{Lemma}\label{first-lemma}
   Let $G$ be a graph and let $v$ be a vertex in $G$ of degree at most five such that there is no $v$-immersion of $K_6$ into $G$. If $G-v$ is $5$-colorable, then so is $G$.
\end{Lemma}

\begin{proof}
  Suppose that $G$ is not $5$-colorable, and let $c$ be a $5$-coloring of $G-v$.
  Then $c$ must assign all five colors to the neighbors of
  $v$ and hence $\deg(v)=5$; since otherwise we can extend $c$ to $G$.
  Let the neighbors of $v$ be $v_1$, $v_2$, $v_3$, $v_4$ and $v_5$; and denote $c(v_i)=i$ for each $i\in\{1,2,3,4,5\}$.

  For each pair of distinct $i$ and $j$ in $\{1,2,3,4,5\}$, let $G_{\{i,j\}}$ denote
  the subgraph of $G-v$ whose vertices are colored by $c$ with $i$ or $j$.
  If, for one such pair of $i$ and $j$, the graph $G_{\{i,j\}}$ has $v_i$ and $v_j$
  in distinct components, then the colors $i$ and $j$ can be switched in one of the
  components so that two neighbors of $v$ are colored the same.
  In this case, the coloring $c$ can be extended to $v$ so that $G$ is
  $5$-colored; a contradiction.
  
  Hence, for each pair of distinct $i$ and $j$, the graph $G-v$ has a path joining $v_i$
  and $v_j$ whose vertices are alternately colored $i$ and $j$ by $c$, and thus $G$ contains
  a $v$-immersion of $K_6$; again, a contradiction.
\end{proof}

\begin{Corollary}[Generalized Five Color Theorem]
        Every graph with crossing number at most two is $5$-colorable.
\end{Corollary}

\begin{proof}
        Suppose not and consider a counterexample $G$ on the minimum number of vertices.
  Proposition~\ref{num_crossings} implies that $m \leq 3n-4$, and so $G$ has a vertex $v$ whose degree is at most five.
  From Corollaries~\ref{k6} and~\ref{onto-greater} we conclude that there is no essential immersion, and hence no $v$-immersion, of $K_6$ into~$G$.
  The minimality of $G$ implies that $G-v$ is $5$-colorable, from which Lemma~\ref{first-lemma} provides the required contradiction.
\end{proof}

Lemma~\ref{first-lemma} establishes that a graph $G$ with $\nu(G) \le 3$ is $5$-colorable if there is no $v$-immersion of $K_6$ into $G$.
The next lemma addresses the case of graphs with $\nu(G)\le3$ for which there is a $v$-immersion of $K_6$ into $G$ for some vertex $v$ in $G$.
The following corollary of a result of Kleitman~\cite{Kleitman1976} will be used in its proof.

\begin{Proposition}\label{k5}
  Every good drawing of $K_5$ has odd number of crossings.
\end{Proposition}

\begin{Lemma}\label{v-immersion}
  If $G$ is a drawing with exactly three crossings and
  $\alpha$ is a $v$-immersion of $K_6$ into $G$ for some vertex $v$ in $G$,
  then $v$ is incident with exactly two crossed edges.
\end{Lemma}

\begin{proof}
  Let $H$ be the subgraph of $G$ that is the image of $K_6$ under $\alpha$, and let $u$ be the vertex in $K_6$ such that $\alpha(u) = v$.
  From Lemma~1.1 and Corollaries~1.5 and~2.2, it follows that $H$ is a good
drawing containing all three crossings of $G$.

  If $v$ were incident with one or three crossed edges in $H$, then $H-v$ would be a good drawing with zero or two crossings with $K_5$ essentially immersed onto it.
  This, together with Lemma~\ref{immersion-crossing}, would imply that there is a good drawing of $K_5$
  with zero or two crossing, which would contradict Proposition~\ref{k5}.
  
  Moreover, if $v$ were incident with no crossed edges in $H$, then $H-v$ would be a drawing with a region $R$ that is incident with all vertices in the set $S=\{\alpha(w)\colon w\in V(K_6-u)\}$.
  The boundary of $R$ then induces a cyclic order on the set $S$, and hence also on
  $V(K_6-u)$.
  If $e$ and $f$ are distinct non-adjacent edges of $K_6-u$ and each joins
  a pair of non-consecutive vertices, then $\alpha(e)$ and $\alpha(f)$ must cross.
  It follows that $H$ would have at least five crossings; a contradiction.
\end{proof}

%%%%%%%%%%%%%%%%%%%%%%%%%%%%%%%%%%%%%%%%%%%%%%%%%%%%%%%%%%%%%%%%%%%%%%%%%%%%%%%%%%%%%%%%

\section{Colorings and Crossings}

Lemmas~\ref{first-lemma} and \ref{v-immersion}, respectively, characterize a graph $G$ when it does not and does contain a $v$-immersion of $K_6$. 
With these, we now proceed to the main theorem.
We will use $\omega(G)$ to denote the {\em clique number} of $G$, that is, the largest $n$ for which $K_n$ is a subgraph of $G$.

\begin{Main}
\label{main-theorem}
  If $\nu(G)\le3$ and $\omega(G)\le5$, then $G$ is $5$-colorable.
\end{Main}

\begin{proof}
  Let $\mathcal{G}$ denote the class of all graphs with crossing number at most three that are not $5$-colorable, and let $G$ be a member of $\mathcal{G}$ with the minimum number of vertices.
  Suppose that $\omega(G) \le 5$ and that $G$ is drawn optimally in the plane.

  If $G$ contains a vertex $v$ of degree less than five, then $G$ is not a minimal
  member of $\mathcal{G}$, since a $5$-coloring of $G-v$ extends to a $5$-coloring
  of $G$.
  Hence, the minimum degree of $G$ is five.
  By Proposition~\ref{num_crossings}, the graph $G$ has at most $3n-3$ edges, and thus
  has at least six vertices of degree five.

  Let $v$ be a vertex of degree five.
  Lemma~\ref{first-lemma} implies that there is a $v$-immersion of $K_6$ into $G$, and Corollary~\ref{onto-greater} implies that the image of $K_6$ in $G$ contains three crossed edges.
  Then Lemma~\ref{v-immersion} implies that two crossed edges of $G$ are incident with $v$. Since $G$ is not $K_6$, it contains a vertex $w$ of degree five not adjacent to $v$. 
  However, Lemma~\ref{first-lemma} implies that there is also a $w$-immersion of $K_6$ into $G$, and so $w$ is also incident with two crossed edges. 
  Since $v$ and $w$ are not adjacent, these two crossed edges are different from the crossed edges incident with $v$, which implies that $G$ contains four crossings; a contradiction.
\end{proof}

We also show that when Theorem~3.1 is applied to a $4$-connected graph $G$ other than $K_6$, then the assumption $\omega(G) \leq 5$ may be discarded.
More precisely, we have:

\begin{Corollary}
If $G$ is $4$-connected, $\nu(G) \le 3$ and $G \neq K_6$, then $G$ is $5$-colorable.
\end{Corollary}
\begin{proof}
Let $G$ be a drawing with at most three crossings of a $4$-connected graph not isomorphic to $K_6$.
We show that $\omega(G) \leq 5$, from which the conclusion follows immediately by Theorem~\ref{main-theorem}.

Suppose, to the contrary, that $G$ has a complete subgraph $K$ on six vertices.
Let $v$ be a vertex of $G$ that is not in $K$, and let $K'$ be the plane drawing
obtained from $K$ by replacing each crossing with a new vertex.
By Corollary~\ref{k6}, all three crossings of $G$ are in $K$, and so $|V(K')|=9$ and $|E(K')|=21$.
Thus $K'$ is a triangulation and so every region of $K$ contains at most three vertices
in its boundary.
But this is impossible, as $G$, being $4$-connected, has four paths from $v$ to vertices
of $K$, with each pair of paths having only $v$ in common.
\end{proof}

Lastly, note that $C_3 \vee C_5$, the graph in which every vertex of $C_3$ is adjacent to every vertex of $C_5$, contains no $K_6$ subgraph and is not 5-colorable.

\begin{figure}[ht]
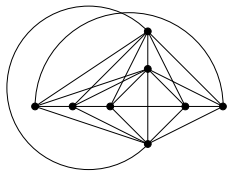

\label{gr2}
\begin{graph}(4.5,2.25)(-.5,1.5)
  \graphnodesize{.1}
  \graphlinewidth{.01}
  \roundnode{1}(.5,2)
  \roundnode{3}(1,2)
  \roundnode{5}(1.5,2)
  \roundnode{7}(2.5,2)
  \roundnode{9}(3,2)
  \roundnode{6}(2,1.5)
  \roundnode{8}(2,2.5)
  \roundnode{10}(2,3)
  
  \edge{1}{6}
  \edge{1}{8}
  \edge{1}{10}
  \edge{3}{6}
  \edge{3}{8}
  \edge{3}{10}
  \edge{5}{6}[\graphlinecolour(0,0,0)]
  \edge{5}{8}[\graphlinecolour(0,0,0)]
  \edge{5}{10}
  \edge{7}{6}[\graphlinecolour(0,0,0)]
  \edge{7}{8}[\graphlinecolour(0,0,0)]
  \edge{7}{10}
  \edge{9}{6}
  \edge{9}{8}
  \edge{9}{10}
  
  \edge{1}{3}
  \edge{3}{5}
  \edge{5}{7}
  \edge{7}{9}
  \bow{1}{9}{.5}
  \edge{6}{8}
  \edge{8}{10}
  \bow{6}{10}{1.25}
\end{graph}
\caption{$C_3 \vee C_5$ drawn with the minimum number of crossings}
\end{figure}

\begin{Proposition}
  The crossing number of $C_3 \vee C_5$ is six.
\end{Proposition}

\begin{proof}
        Let $G$ be an optimal drawing of $K \vee L$, where $K$ and $L$ are
        cycles on, respectively, three and five vertices.
        Suppose that $G$ has fewer than six crossings.
        Note that $G\backslash (E(K)\cup E(L))$ is isomorphic to $K_{3,5}$,
        which has crossing number four \cite{Kleitman1970}.
        This implies that the edges of $K\cup L$ are involved in at most
        one crossing, and thus $L$ has at most three regions, one of which
        contains $K$.
        Thus at least one region of $L$ avoids $K$ and has two non-adjacent
        vertices of $L$ in its boundary.
        These two vertices of $L$ can be joined by a new edge that crosses no
        edges of $G$ thereby creating a graph with $8$ vertices, $24$ edges, and
        $5$ crossings; a contradiction to Proposition~\ref{num_crossings}.
        Hence, $G$ has six crossings.
        Figure 2 shows a drawing which achieves this bound, proving that $\nu(C_3 \vee C_5) = 6$.
\end{proof}

We do not currently know whether the Main Theorem~\ref{main-theorem} extends to graphs with four or five crossings, and hence conclude with the following question:

\begin{Question}
Does a graph $G$ have a $5$-coloring if $\nu(G) \le 5$ and $\omega(G) \le 5$?
\end{Question}

\bibliographystyle{amsalpha}

\end{document}